\documentclass[12pt,twoside]{article}
\usepackage{amsmath,amsbsy,amsfonts}

\newtheorem{theorem}{Theorem}[section]
\newtheorem{lemma}[theorem]{Lemma}

\newcommand{\qed}{\hfill\rule{2mm}{2mm}}

\title{Existence of standing waves solution for a Nonlinear Schr\"{o}dinger equations in $\mathbb{R}^{N}$}

\author{Claudianor O. Alves\footnote{Research of C. O. Alves partially supported by  CNPq 304036/2013-7  and INCT-MAT, e-mail:coalves@dme.ufcg.edu.br}  \\
Universidade Federal de Campina Grande\\
Unidade Acad\^emica de Matem\'atica\\
CEP:58429-900, Campina Grande - PB, Brazil.}

\date{}

\begin{document}
\maketitle



\begin{abstract}

In this paper, we investigate the existence of positive solution for
the following class of elliptic equation
$$
- \epsilon^{2}\Delta u +V(x)u= f(u) \,\,\,\, \mbox{in} \,\,\, \mathbb{R}^{N},
$$
where $\epsilon >0$ is a positive parameter, $f$ has a subcritical growth and $V$ is a positive potential verifying some conditions.

\end{abstract}


\section{Introduction}

In recent years, many authors have considered the existence of solution for the following class of elliptic equation 
$$
\left\{\begin{array}{l}
- \epsilon^{2} \Delta u +V(x)u= f(u) \,\,\,  \mbox{ in } \mathbb R^N, \\
u>0 \;\;\; \mbox{in} \;\;\; \mathbb{R}^{N},\\
 u \in  H^{1}(\mathbb R^N),
\end{array}\right. \eqno (P)_\epsilon
$$
where $\epsilon >0$ is a positive parameter, $V:\mathbb R^N \to \mathbb R$ and $f:\mathbb R \to \mathbb R$ are continuous functions with $V$ being a nonnegative function and $f$ having a subcritical or critical growth. The existence and concentration of
positive solutions for general semilinear elliptic equations $(P)_\epsilon$ for the case $N \geq 2$ have been extensively
studied, see for example, Ackermann and Szulkin \cite{Ackermann}, Alves, do \'O and Souto \cite{AOS}, Bartsch, Pankov and Wang \cite{BPW}, do \'O and Souto \cite{OS}, del Pino and Felmer \cite{Pino, PFM2}, del Pino, Felmer and Miyagaki \cite{PFM2}, Floer and Weinstein \cite{FW}, Oh \cite{Oh1}, Rabinowitz \cite{r}, Wang \cite{W} and their references.

Knowledge of the solutions of $(P)_\epsilon$ has a great importance for studying the existence of {\it standing wave solutions} for nonlinear Schr\"{o}dinger equation
$$
i \epsilon \displaystyle \frac{\partial \Psi}{\partial t}=- \epsilon^{2}\Delta
\Psi+W(z)\Psi-f(\Psi)\,\,\, \mbox{for all}\,\,\, z \in \mathbb{R}^{N}, \eqno{(NLS)}
$$
which are solutions of the form $\Psi(x,t)=\exp(-iEt/\epsilon)u(x)$, where $u$ is a solution of $(P)_\epsilon$. The equation $(NLS)$ is one of the main objects of the quantum physics, because it appears in problems involving nonlinear optics, plasma physics and condensed matter physics, see \cite{FW} and \cite{Oh1} for more details about this subject.

In a seminal paper, Rabinowitz \cite{r} introduced the following condition on $V$ 
$$
0<\inf_{z \in \mathbb{R}^{N}}V(z)< \liminf_{|z| \to +\infty}V(z). \eqno{(V_8)}
$$
Later Wang \cite{W} showed that these solutions concentrate at global minimum points of $V$ as  $\epsilon$ tends to 0.

In \cite{Pino}, del Pino and Felmer established  the existence of positive  solutions which concentrate around local minimum of $V$, by introducing a penalization method. More precisely, they assumed that there is an open and bounded set $\mathcal{O}$ compactly contained  in $\mathbb{R}^{N}$ such that
$$
0< \gamma \leq V_0 =\inf_{z\in \mathcal{O}}V(z)< \min_{z \in
	\partial \mathcal{O}}V(z). \eqno{(V_{1})}
$$
Motivated by this result, Alves, do \'O and Souto \cite{AOS} and do \'O and Souto \cite{OS} studied the same type of problem with $f$ having critical growth for $N \geq 3$ and exponential critical growth for $N=2$ respectively.

In \cite{PFM}, del Pino, Felmer and Miyagaki considered the case where potential $V$ has a geometry like saddle, essentially they assumed the following conditions on $V$: First of all, they fixed two subspaces $X,Y  \subset \mathbb{R}^{N}$ such that
$$
\mathbb{R}^{N}=X \oplus Y.
$$  
By supposing that $V$ is bounded,  they fixed $c_0,c_1>0$  satisfying \\
$$
\displaystyle c_0=\inf_{z \in \mathbb{R}^{N}}V(z)>0
$$
and
$$
c_1=\displaystyle \sup_{x \in X}V(x).
$$
Furthermore, they also supposed that $V \in C^{2}(\mathbb{R}^{N})$ and it  verifies the following geometric conditions:  \\

\noindent $(V_1)$  
$$
c_0=\inf_{R>0}\sup_{x \in \partial B_R(0) \cap X }V(x)<\inf_{y \in Y}V(y).
$$
\noindent $(V_2)$ \quad The functions $V, \frac{\partial V}{\partial x_i}$ and $\frac{\partial^{2} V}{\partial x_i \partial x_j}$ are bounded in $\mathbb{R}^{N}$ for all $i,j \in \{1,...,N\}$.  \\

\noindent $(V_3)$ \quad $V$ satisfies the Palais-Smale condition, that is, if $(x_n) \subset \mathbb{R}^{N}$ is a sequence such that $(V(x_n))$ is bounded and $\nabla V(x_n) \to 0$, then $(x_n)$ possesses a convergent subsequence in $\mathbb{R}^{N}$. \\  

Using the above conditions on $V$, and supposing that 
$$
c_1<2^{\frac{2(p-1)}{N+2-p(N-2)}}c_0,
$$
del Pino, Felmer and Miyagaki showed  the existence of positive solutions for the following problem  
$$
- \epsilon^{2} \Delta{u} + V(z)u=|u|^{p-2}u
\ \ \mbox{in} \ \ \mathbb{R}^{N},
$$
where $p \in (2,2^{*})$ if $N \geq 3$ and $p \in (2,+\infty)$ if $N=1,2$, for $\epsilon>0$ small enough. The main tool used was the variational method, more precisely, the authors found critical points of the functional
$$
E_{\epsilon}(u)=\int_{\mathbb{R}^{N}}(\epsilon^{2}|\nabla u|^{2}+|u|^{2})\,dx
$$ 
on the manifold
$$
\mathcal{M} =\left\{u \in H^{1}(\mathbb{R}^{N}) \cap P \,:\,\int_{\mathbb{R}^{N}}|u|^{p}\,dx=1 \right\},
$$
where $P$ denotes the cone of nonnegative functions of $H^{1}(\mathbb{R}^{N})$.  Recently, in \cite{AlvesNovo}, Alves has studied the same type of problem with $f$ having an exponential critical growth and $N=2$. 

In \cite{PFM2}, del Pino and Felmer have considered the following assumptions on $V$: \\

\noindent $(V_4)$  $V$ is of class $C^{1}$ and there is $\alpha>0$ such that
$$
V(z) \geq \alpha, \quad \forall z \in \mathbb{R}^{N}.
$$

Locally, it was fixed  an open and bounded set $D \subset \mathbb{R}^{N}$ and subsets $B_0,B \subset D$ with $B$ connected. Using these sets, we denoted by $\Gamma$ the class of all continuous function $\phi:B \to D$ with the property that
$\phi(y)=y$ for all $y \in B_0$. Define the min-max value $c$ as
\begin{equation} \label{Nivel}
c=\inf_{\phi \in \Gamma}\sup_{y \in B}V(\phi(y)),
\end{equation}
and assume additionally

\noindent $(V_5)$
$$
\sup_{y \in B_0}V(y)<c.
$$

\noindent $(V_6)$ \,\, For all $\phi \in \Gamma$, $\phi(B) \cap \{y \in D\,:\,V(y) \geq c \} \not= \emptyset$. \\

\noindent $(V_7)$ For all $y \in \partial D$ such that $V(y)=c$, one has $\partial_\nu V(y) \not= 0$, where $\partial_\nu$ denotes the tangential derivative.

\vspace{0.5 cm}

Motivated by the above papers, in the present article  we show the existence of solution for $(P)_\epsilon$, by considering two new classes of potential $V$, namely : \\

\noindent {\bf Classe 1: Potential $V$ verifies the Palais-Smale condition.}
\begin{enumerate}
\item[$(A_{0})$\ ]  There is $V_0>0$ such that $V(x) \geq V_0,$ \,\,\, $\forall x \in \mathbb R^N$.
\end{enumerate}
\begin{enumerate}
\item[$(A_{1})$\ ] $V \in C^{2}(\mathbb{R}^{N})$ and $V, \frac{\partial V}{\partial x_i}$ and $\frac{\partial^{2} V}{\partial x_i \partial x_j}$ are bounded in $\mathbb{R}^{N}$ for all $i,j \in \{1,...,N\}$. 
\end{enumerate}
\begin{enumerate}
\item[$(A_{2})$\ ] $V$ verifies the Palais-Smale condition, that is,  if $(x_n) \subset \mathbb{R}^{N}$ is such that $(V(x_n))$ is bounded and $\nabla V(x_n) \to 0$, then $(x_n)$ possesses a convergent subsequence in $\mathbb{R}^{N}$. 
\end{enumerate}

\noindent {\bf Classe 2: Potential $V$ does not have critical point on the boundary of some bounded domain. } \\

In this class of potential, we suppose that $V$ verifies $(A_0)-(A_1)$ and the following additional condition: 
\begin{enumerate}
\item[$(A_{3})$\ ]  There is a bounded domain $\Lambda \subset \mathbb{R}^{N}$, such that 
$ \nabla V(x) \not=0$ for all $ x \in \partial \Lambda. $
\end{enumerate}

Related to function $f$, we assume that 
\begin{enumerate}
\item[$(f_1)$\ ]
$\displaystyle\limsup_{s\rightarrow 0 ^+} \frac{f(s)}{s}=0.$

\item[$(f_2)$\ ]  There exists $p \in (2,2^*)$, such that
$$
\limsup_{s\rightarrow +\infty} \frac{f(s)}{s^{p-1}} =
0;
$$

\item[$(f_3)$\ ] There exists $\theta>2$  such that
$$
0<\theta F(s)\leq sf(s) \,\,\,\,\,\,\,\,\,  \forall s >0; 
$$
\end{enumerate}

The statement of our main result is the following    

\begin{theorem}{\label{T1}}
Suppose that $V$ belongs to Classe 1 or 2 and $f$  satisfies $(f_1)-
(f_3)$. Then, problem $(P)_\epsilon$ has a positive solution for $\epsilon>0$ small enough. 
\end{theorem}

\vspace{1mm}

In the proof of Theorem \ref{T1}, we will use variational methods, more precisely, the Mountain Pass Theorem due to Ambrosetti and Rabinowitz \cite{ARMPT} combined with some arguments developed by del Pino and Felmer \cite{Pino}, for more details see Section 2. 

The paper is organized as follows. In the next section, inspired by \cite{Pino}, we study the existence of solution for a class of auxiliary problem. In Section 3, we prove the main theorem supposing that $V$ belongs to Class 1, while in Section 4, we consider the case that $V$ belongs to Class 2. Finally, in Section 5, we make some final considerations for elliptic problems with $f$ having critical growth for $N \geq 3$ and exponential critical growth for $N=2$.

\section{del Pino and Felmer's approach} 

In this section, following an idea  found in del Pino and Felmer \cite{Pino}, we will study the existence of solution for a special class of elliptic problem associated with $(P)_\epsilon$. 

Since we intend to prove the existence of positive solutions, hereafter  we
consider 
$$
f(s)=0, \quad \forall s \leq 0.
$$

Using the change variable $v(x)=u(\epsilon x)$, it is possible to prove that $(P)_\epsilon$ is equivalent to the following problem
$$
\left\{\begin{array}{l}
-\Delta u +V(\epsilon x)u= f(u) \,\,\,  \mbox{ in } \mathbb R^N, \\
u>0 \;\;\; \mbox{in} \;\;\; \mathbb{R}^{N},\\
 u \in  H^{1}(\mathbb R^N).
\end{array}\right. \eqno (P)'_\epsilon
$$
Therefore, in what follows, we prove the existence of positive solution for $(P)'_\epsilon$. To this end, we start observing that from $(A_{0})$, we can work in $H^{1}(\mathbb{R}^{N})$ with the norm 
$$
\|u\|=\left(\int_{\mathbb R^N}(|\nabla u|^2+V(\epsilon x)|u|^2)dx\right)^{\frac{1}{2}},
$$
which is equivalent to usual norm. 

The Euler-Lagrange functional associated with $(P)'_\epsilon$ is given by
$$
I_\epsilon(u)=\frac 12 \|u\|^2 - \int_{\mathbb R^N}F(u)dx, \;\; \forall u \in H^{1}(\mathbb{R}^{N}).
$$
From the conditions on $f$, the functional $I_\epsilon \in C^1(H^{1}(\mathbb{R}^{N}),\mathbb R)$
and its Gateaux derivative is
$$
I'_\epsilon(u)v= \int_{\mathbb R^N} (\nabla u \nabla v + V(\epsilon x) uv)dx -
\int_{\mathbb R^N}f(u)vdx, \,\, \forall u,v \in H^{1}(\mathbb{R}^{N}).
$$
It is easy to check that the critical points of $I_\epsilon$ are weak solutions of $(P)'_\epsilon$.

In the sequel, let us denote by $I_{\infty}:H^{1}(\mathbb{R}^{N}) \to \mathbb{R}$ the functional  
$$
I_\infty(u)=\frac 12 \int_{\mathbb{R}^{N}}(|\nabla u|^2+V_\infty |u|^2)dx -
\int_{\mathbb{R}^{N}}F(u)dx,
$$
where 
$$
V_\infty=\max_{x \in \mathbb{R}^{N}}V(x).
$$ 
Furthermore, let us denote by $c_\infty$ the mountain level associated with $I_{\infty}$, that is,
$$
c_\infty= \inf_{\gamma \in \Gamma}\max_{t\in [0,1]} I_\infty(\gamma(t))
$$
where 
\begin{equation} \label{XQ}
\Gamma = \{ \gamma \in C([0,1], H^{1}(\mathbb{R}^{N})): \gamma (0)=0 \,\,\, \mbox{and} \,\,\, 
I_\infty(\gamma(1))<0 \}.
\end{equation}

Here, we would like point out that $c_\infty$ depends only on $V_\infty$, $\theta$
and $f$.

\subsection{An auxiliary problem}

Given a bounded domain $\Omega \subset \mathbb{R}^{N}$, we fix the numbers $k=\frac {2\theta}{\theta-2}>2$ and  $a>0$ be the value at which 
$$
\frac{f(a)}{a}=\frac{V_0}{k},
$$ 
where $V_0>0$ was given in $(A_0)$. Using these numbers, we set the functions
$$
\tilde{f}(s)=
\left\{
\begin{array}{l}
0, \quad s \leq 0, \\
f(s), \quad 0 \leq s \leq a, \\
\frac{V_0}{k}s, \quad s \geq a
\end{array}
\right.
$$
and
$$
g(x,s)=\chi(x)f(s)+(1-\chi)\tilde{f}(s), \quad \forall (x,s) \in \mathbb{R}^{N} \times \mathbb{R},
$$
where  $\chi$ denotes the characteristic function associated with $\Omega$, that is, 
$$
\chi(x)=
\left\{
\begin{array}{l}
1, \quad x \in \Omega \\
0, \quad x \in \Omega^c.
\end{array}
\right.
$$

Using the above functions, we will study the existence of positive solution for the following problem
$$
\left\{
\begin{array}{l}
-\Delta{u}+V(\epsilon x)u=g_\epsilon(x,u), \quad x \in \mathbb{R}^{N}, \\
\mbox{}\\
u \in H^{1}(\mathbb{R}^{N}),
\end{array}
\right.
\leqno{(AP)_{\epsilon}}
$$
where 
$$
g_\epsilon(x,s)=g(\epsilon x,s), \quad \forall (x,s) \in \mathbb{R}^{N} \times \mathbb{R}.
$$

The above problem is strongly related with $(P)'_\epsilon$, because if $u$ is a solution of $(AP)_{\epsilon}$ verifying
$$
u(x)<a, \quad \forall x \in \mathbb{R}^{N} \setminus \Omega_\epsilon,
$$
where $\Omega_\epsilon={\Omega}/{\epsilon}$, then $u$ will be a solution for $(P)'_\epsilon$.

Associated with $(AP)_{\epsilon}$, we have the energy functional $J_\epsilon:H^{1}(\mathbb{R}^{N}) \to \mathbb{R}$ given by
$$
J_\epsilon(u)=\frac{1}{2}\int_{\mathbb{R}^{N}}(|\nabla u|^{2}+V(\epsilon x)|u|^{2})\,dx-\int_{\mathbb{R}^{N}}G_\epsilon(x,u)\,dx,
$$
where
$$
G_\epsilon(x,s)=\int_{0}^{s}g_\epsilon(x,t)\,dt, \quad \forall (x,s) \in \mathbb{R}^{N} \times \mathbb{R}.
$$

From the conditions on $f$, and hence on $g$, $J_\epsilon \in C^1(H^{1}(\mathbb{R}^{N}),\mathbb R)$
with
$$
J'_\epsilon(u)v= \int_{\mathbb R^N} (\nabla u \nabla v + V(\epsilon x) uv)dx -
\int_{\mathbb R^N}g_\epsilon(x,u)vdx, \, \forall u,v \in H^{1}(\mathbb{R}^{N}).
$$
Thus, critical points of $J_\epsilon$ correspond to weak solutions of $(AP)_{\epsilon}$.

Repeating the same arguments found in \cite{Pino}, it is easy to see that $J_\epsilon$ verifies the hypotheses of the Mountain Pass Theorem due to Ambrosetti and Rabinowitz \cite{ARMPT} for all $\epsilon >0$. Therefore, there is $u_\epsilon \in H^{1}(\mathbb{R}^{N})$ such that
$$
J_\epsilon(u_\epsilon)=c_\epsilon>0 \quad \mbox{and} \quad J'_\epsilon(u_\epsilon)=0,
$$ 
where
$$
c_\epsilon= \inf_{\gamma \in \Gamma}\max_{t\in [0,1]} J_\epsilon(\gamma(t))
$$
with  
$$
\Gamma = \{ \gamma \in C([0,1],H^{1}(\mathbb{R}^{N})): \gamma (0)=0 \,\,\, \mbox{and} \,\,\, 
J_\epsilon(\gamma(1))<0 \}.
$$
Observing that 
$$
J_\epsilon(u) \leq I_{\infty}(u) \quad \forall u \in H^{1}(\mathbb{R}^{N}),
$$ 
we ensure that 
\begin{equation} \label{minimax}
c_\epsilon \leq c_\infty, \quad \forall \epsilon >0.
\end{equation}

The lemma below establishes an important estimate from above for the $H^{1}$-norm of the family $(u_\epsilon)$.
\begin{lemma}{\label{lm4}}
For all $\epsilon >0$, the solution $u_\epsilon$ of $(AP)_{\epsilon}$ 
satisfies the estimate
\[
\|u_\epsilon\|^2\leq 2kc_\infty.
\]
\end{lemma}

\noindent {\bf{Proof:}}  Using the fact that $u_\epsilon$ is a critical point of $J_\epsilon$, we must have 
$$
c_\epsilon=J_\epsilon(u_\epsilon)=J_\epsilon(u_\epsilon)-\frac 1\theta J'_\epsilon(u_\epsilon)u_\epsilon \geq \frac {(\theta-2)}{4\theta}\|u_\epsilon \|^2 \geq \frac 1{2k}\|u_\epsilon\|^2.
$$
Now, the result follows combining the above inequality with (\ref{minimax}). \qed

\vspace{0.5 cm}

Here, we would like point out that in Lemma \ref{lm4}, the norm $\|u_\epsilon \|$ is bounded from above, by a
constant that depends only of $V_\infty$, $\theta$ and $f$, then the  constant does not depend on $\epsilon >0$.

 \vspace{1mm}



\section{Proof of Theorem \ref{T1}: The Class 1}

In this section, we will prove the Theorem \ref{T1}, by supposing that $V$ belongs to Class 1. To this end, we will use the results obtained in Section 2  fixing 
$$
\Omega = B_{R_\epsilon}(0),
$$ 
where $R_\epsilon=\frac{1}{\epsilon}$. Consequently, we know that there is a solution $u_\epsilon \in H^{1}(\mathbb{R}^{N})$ for  $(AP)_\epsilon$.

In what follows, our goal is to prove that there is $\epsilon_0>0$ such that
$$
u_\epsilon(x) <a, \quad \forall x \in \mathbb{R}^{N} \setminus B_{\frac{R_\epsilon}{\epsilon}}(0) \quad \mbox{and} \quad \forall \epsilon \in (0,\epsilon_0).
$$

\begin{lemma}{\label{lm6}} The function $u_\epsilon$ verifies the following estimate 
$$
\max_{x \in \partial B_{\frac{R_\epsilon}{\epsilon}}(0)}u_\epsilon(x) \to 0, \quad \epsilon \to 0.
$$
\end{lemma}

\noindent {\bf{Proof:}}  Assume by contradiction that there is $\epsilon_n \to 0$ and $\gamma >0$ such that 
$$
\max_{x \in \partial B_{\frac{R_{\epsilon_n}}{\epsilon_n}}(0)}u_{n}(x) \geq \gamma \quad  \forall n \in  \mathbb{N},
$$
where $u_n=u_{{\epsilon_n}}$. From now on, we fix $x_n \in \partial B_{\frac{R_{\epsilon_n}}{\epsilon_n}}(0)$ satisfying
$$
u_n(x_n)=\max_{x \in \partial B_{\frac{R_{\epsilon_n}}{\epsilon_n}}(0)}u_{{\epsilon_n}}(x).
$$
Therefore,
$$
u_n(x_n) \geq \gamma, \quad  \forall n \in  \mathbb{N}.
$$
By Lemma \ref{lm4}, $(u_n)$ is bounded in $H^{1}(\mathbb{R}^{N})$. Thereby, setting $w_n=u_n(\cdot+x_n)$, we can guarantee that $(w_n)$ is also bounded in $H^{1}(\mathbb{R}^{N})$ and it satisfies 
$$
\left\{\begin{array}{l}
- \Delta w_n +V(\epsilon_n x+\epsilon_n x_n)w_n= g(\epsilon_n x+\epsilon_n x_n,w_n), x \in \mathbb R^N, \\
\mbox{}\\
 u \in H^{1}(\mathbb R^N).
\end{array}\right. 
$$
Using bootstrap arguments, it is possible to show that $(w_n)$ converges uniformly on compact set for its weak limit $w \in  H^{1}(\mathbb{R}^{N})$. Then,  $w \in C(\mathbb{R}^{N})$ and $w(0) \geq \gamma$, implying that $w\not=0$.  Moreover, by $(A_1)$, there exists a subsequence of $(\epsilon_n x_n)$, still denote by itself, such that 
$$
\alpha_1=\displaystyle \lim_{n \to +\infty}V(\epsilon_n x_n),
$$ 
for some $\alpha_1>0$. Since for each $\phi \in H^{1}(\mathbb{R}^{N})$, the equality below holds
$$
\int_{\mathbb{R}^{N}}\nabla w_n \nabla \phi \, dx + \int_{\mathbb{R}^{N}}V(\epsilon_n x +\epsilon_n x_n )w_n \phi \, dx - \int_{\mathbb{R}^{N}}g(\epsilon_n x+\epsilon_n x_n,w_n)\phi \, dx = o_n(1)\|\phi\|,
$$
taking the limit of $n \to +\infty$, let us deduce that $w$ is a nontrivial solution of the problem
\begin{equation} \label{equacao}
\Delta{u}-\alpha_1u+\tilde{g}(x,u)=0, \quad x \in \quad \mathbb{R}^{N},
\end{equation}
where 
$$
\tilde{g}(x,s)=\tilde{\chi}(x)f(s)+(1-\tilde{\chi}(x))\tilde{f}(s),
$$
for some $\tilde{\chi} \in L^{\infty}(\mathbb{R}^{N})$. Thus, by regularity theory, $w \in L^{\infty}(\mathbb{R}^{N}) \cap H^{2}(\mathbb{R}^{N})$. 

For each $j \in \mathbb{N}$, there is $\phi_j \in C^{\infty}_{0}(\mathbb{R}^{N})$ such that
$$
\|\phi_j -w\| \leq {1}/{j},
$$
that is,
$$
\|\phi_j -w\|=o_j(1).
$$

Using $\frac{\partial \phi_j}{\partial x_i}$ as a test function, we get
$$
\int_{\mathbb{R}^{N}}\nabla w_n \nabla \frac{\partial \phi_j}{\partial x_i}\,dx+ \int_{\mathbb{R}^{N}}V(\epsilon_n x+\epsilon_n x_n)w_n\frac{\partial \phi_j}{\partial x_i}\,dx-
\int_{\mathbb{R}^{N}}g(\epsilon_n x+\epsilon_n x_n,w_n)\frac{\partial \phi_j}{\partial x_i}\,dx=o_n(1).
$$
Now, using well known arguments,
$$
\int_{\mathbb{R}^{N}}\nabla w_n \nabla \frac{\partial \phi_j}{\partial x_i}\,dx=\int_{\mathbb{R}^{N}}\nabla w \nabla \frac{\partial \phi_j}{\partial x_i}\,dx+o_n(1)
$$
and
$$
\int_{\mathbb{R}^{N}}g(\epsilon_n x+\epsilon_n x_n,w_n)\frac{\partial \phi_j}{\partial x_i}\,dx=\int_{\mathbb{R}^{N}}\tilde{g}(x,w)\frac{\partial \phi_j}{\partial x_i}\,dx+o_n(1).
$$
Gathering the above limit with (\ref{equacao}), we find  
$$
\limsup_{n \to +\infty}\left|\int_{\mathbb{R}^{N}}(V(\epsilon_n x+\epsilon_n x_n)-V(\epsilon_n x_n))w_n\frac{\partial \phi_j }{\partial x_i}\,dx\right|=0.
$$
As $\phi_j$ has compact support, the above limit gives 
$$
\limsup_{n \to +\infty}\left|\int_{\mathbb{R}^{N}}(V(\epsilon_n x+\epsilon_n x_n)-V(\epsilon_n x_n))w\frac{\partial \phi_j}{\partial x_i}\,dx\right|=0.
$$
Now, recalling that $\frac{\partial w}{\partial x_i} \in L^{2}(\mathbb{R}^{N})$, we have that $(\frac{\partial \phi_j}{\partial x_i})$ is bounded in $L^{2}(\mathbb{R}^{N}).$ Hence, 
$$
\limsup_{n \to +\infty}\left|\int_{\mathbb{R}^{N}}(V(\epsilon_n x+\epsilon_n x_n)-V(\epsilon_n x_n))\phi_j\frac{\partial \phi_j}{\partial x_i}\,dx\right|=o_j(1),
$$
and so,
$$
\limsup_{n \to +\infty}\left|\frac{1}{2}\int_{\mathbb{R}^{N}}(V(\epsilon_n x+ \epsilon_n x_n)-V(\epsilon_n x_n))\frac{\partial (\phi_j^{2})}{\partial x_i}\,dx\right|=o_j(1).
$$
Using Green's Theorem together with the fact that $\phi_j$ has compact support, we obtain the limit below
$$
\limsup_{n \to +\infty}\left|\int_{\mathbb{R}^{N}}\frac{\partial V}{\partial x_i}(\epsilon_n x+ \epsilon_n x_n) \, \phi_j^{2}\right|\,dx=o_j(1),
$$
which combined with $(A_1)$ loads to
$$
\limsup_{n \to +\infty}\left|\frac{\partial V}{\partial x_i}(\epsilon_n x_n)\int_{\mathbb{R}^{N}}|\phi_j|^{2}\,dx\right|=o_j(1).
$$
As 
$$
\int_{\mathbb{R}^{2}}|\phi_j|^{2}\,dx \to \int_{\mathbb{R}^{N}}|w|^{2}\,dx >0\quad \mbox{as} \quad j \to +\infty,
$$
it follows that
$$
\limsup_{n \to +\infty}\left|\frac{\partial V}{\partial x_i}(\epsilon_n x_n)\right|=o_j(1), \quad \forall i \in \{1,....,N\}.
$$
Since $j$ is arbitrary, we derive that 
$$
\nabla V(\epsilon_n x_n) \to 0 \quad \mbox{as} \quad n \to \infty.
$$
Therefore,  $(\epsilon_n x_n)$ is a $(PS)_{\alpha_1}$ sequence for $V$, which is an absurd, since by $(A_2)$,  $V$ satisfies the $(PS)$ condition and $(\epsilon_n x_n)$ does not have any convergent subsequence in $\mathbb{R}^{N}$, because
$$
|\epsilon_n x_{\epsilon_n}|=R_{\epsilon_n}=\frac{1}{\epsilon_n} \to +\infty, \quad \mbox{as} \quad n \to +\infty.
$$
\qed

\noindent {\bf{Proof of Theorem \ref{T1} ( conclusion ):}} From Lemma \ref{lm6}, there is $\epsilon_0>0$ such that
$$
\max_{x \in \partial B_{\frac{R_\epsilon}{\epsilon}}(0)}u_\epsilon(x)<a, \quad \forall \epsilon \in (0, \epsilon_0).
$$
Considering the function 
$$
\tilde{u}_\epsilon(x)=
\left\{
\begin{array}{l}
0, \quad x \in \overline{B}_{\frac{R_\epsilon}{\epsilon}}(0), \\
\mbox{}\\
(u_\epsilon-a)^{+}(x), \quad x \in \mathbb{R}^{N} \setminus {B_{\frac{R_\epsilon}{\epsilon}}(0)},
\end{array}
\right.
$$
it follows that $\tilde{u}_\epsilon \in H^{1}(\mathbb{R}^{N})$. Thereby, $J'_\epsilon(u_\epsilon)\tilde{u}_\epsilon=0$, or equivalently,
$$
\int_{\mathbb{R}^{N}}\nabla {u}_\epsilon \nabla \tilde{u}_\epsilon \, dx +\int_{\mathbb{R}^{N}}V(\epsilon x){u}_\epsilon \tilde{u}_\epsilon \, dx = \int_{\mathbb{R}^{N}}g_\epsilon(x,{u}_\epsilon)\tilde{u}_\epsilon \, dx.
$$
Now, using the definition of $g_\epsilon$, it is possible to prove that $\tilde{u}_\epsilon \equiv 0$. From this, 
$$
{u}_\epsilon(x) \leq a, \quad \forall x \in \mathbb{R}^{N} \setminus {B_{\frac{R_\epsilon}{\epsilon}}(0)},
$$
showing that ${u}_\epsilon$ is a solution for $(P)'_\epsilon $. \qed

\section{Proof of Theorem \ref{T1}: The Class 2}

In this section, we will prove the Theorem \ref{T1} for the case that $V$ belongs to Class 2. However, we will use the results showed in Section 2 with 
$$
\Omega = \Lambda.
$$ 
Then, we also have a solution $u_\epsilon \in H^{1}(\mathbb{R}^{N})$ for  $(AP)_\epsilon$.

Next, we will show that there is $\epsilon_0>0$ such that
$$
u_\epsilon(x) <a, \quad \forall x \in \mathbb{R}^{N} \setminus \Lambda_\epsilon \quad \mbox{and} \quad \forall \epsilon \in (0,\epsilon_0).
$$

\begin{lemma}{\label{lm60}} The function $u_\epsilon$ verifies the following estimate 
$$
\max_{x \in \partial \Lambda_\epsilon}u_\epsilon(x) \to 0, \quad \epsilon \to 0.
$$
\end{lemma}

\noindent {\bf{Proof:}}  Using the same type of arguments explored in the proof of Lemma \ref{lm6}, we find a sequence $(x_n) \ \subset \partial \Lambda_{\epsilon_n} $, with $\epsilon_n \to 0$, satisfying
$$
\nabla V(\epsilon_n x_n) \to 0.
$$
Since $(\epsilon_n x_n) \subset \partial \Lambda$, and $\partial \Lambda$ is a compact set in $\mathbb{R}^{N}$, we can assume that there is $x_0 \in \partial \Lambda$ such that
$$
\epsilon_n x_n \to x_0 \quad \mbox{in} \quad \mathbb{R}^{N}.
$$
Gathering the above limits with the fact that $V \in C^{1}(\mathbb{R}^{N},\mathbb{R})$, we get 
$$
x_0 \in \partial \Lambda \quad \mbox{and} \quad \nabla V(x_0)=0,
$$
contradicting $(A_3)$.  \qed

\vspace{0.5 cm}

\noindent {\bf{Proof of Theorem \ref{T1} ( conclusion ):}}  The conclusion of the proof follows as in Section 3.

\section{Final considerations}

In this section, we would like point out that the arguments explored in the present paper can be applied to study the existence of solution for elliptic problems with critical growth for $N \geq 3$ and exponential critical growth for $N=2$, with the following hypotheses: \\

\noindent {\bf Critical growth for $N\geq 3 $:}

$$
\left\{
\begin{array}{l}
-\epsilon^{2}\Delta{u}+V(x)u=\lambda |u|^{q-2}u+|u|^{2^{*}-2}u\,\,\, \mbox{in} \,\,\, \mathbb{R}^{N},  \\
u \in H^{1}(\mathbb{R}^{N}),
\end{array}
\right.
\eqno{(P_\epsilon)_*}
$$
where $\epsilon, \lambda >0$ are positive parameters, $q \in (2,2^{*})$ and $2^{*}=\frac{2N}{N-2}$. \\

The main result associated with this class of problem is the following

\begin{theorem} \label{T2} Assume that $V$ belongs to Class 1 or 2. Then,  there is $\epsilon_0>0$ such that \\
\noindent a) \, If $N\geq 4$,  $(P_\epsilon)_*$ has a positive solution for all $\epsilon \in (0,\epsilon_0]$ and $\lambda >0$. \\
\noindent b) \, If $N=3$, there is $\lambda^{*}>0$, which is independent of $\epsilon_0>0$, such that  $(P_\epsilon)_*$ has a positive solution for all $\epsilon \in (0,\epsilon_0]$ and  $\lambda \geq \lambda^{*}$. 
\end{theorem}

In the proof of this result, we adapt the arguments found in \cite{AOS}, because in that paper, problem $(P_\epsilon)_*$ has been considered with the same condition on $V$ considered in \cite{Pino}. \\

\noindent {\bf Critical growth for $N=2 $:} \\

Hereafter, $f \in C^{1}(\mathbb{R})$ and it satisfies the following conditions:

\begin{enumerate}
\item[$(f_1)$] There is $C>0$ such that
$$
|f(s)| \leq Ce^{4\pi |s|^2}\ \ \mbox{for all}\ \ s\in \mathbb{R}.
$$
\item[$(f_2)$] $\displaystyle \lim_{s\rightarrow
0}\dfrac{f(s)}{s}=0$.
\item[$(f_3)$] There is $\theta>2$ such that
$$
0<\theta F(s):=\theta\int_0^{s}f(t)dt\leq sf(s),\ \ \mbox{for all}\ \ s \in \mathbb{R} \setminus \{0\}. 
$$
\item[$(f_4)$] There exist constants $p>2$ and $C_p>0$ such  that
$$
f(s)\geq C_p s^{p-1}\ \ \mbox{for all}\ \ s >0,
$$
where
$$
C_p>\left[\beta_p\left(\frac{2\theta}{\theta-2}\right)\frac{1}{\min\{1,c_0\}}\right]^{(p-2)/2},
$$
with
$$
\beta_p=\inf_{\mathcal{N}_\infty}J_\infty,
$$
$$
\mathcal{N}_\infty=\{u\in H^1(\mathbb{R}^{2}) \setminus \{0\}:\ J_{\infty}'(u)u=0\}
$$
and
$$
J_{\infty}(u)=\dfrac{1}{2}\int_{\mathbb{R}^{2}}\left(|\nabla u|^2+|V|_\infty|u|^2\right)dx-\dfrac{1}{p}\int_{\mathbb{R}^{2}}|u|^pdx .
$$
\end{enumerate}

The main result related to above hypotheses is the following

\begin{theorem} \label{T3} Assume that $V$ belongs to Class 1 or 2 and $f$ verifies $(f_1)-(f_4)$. Then, problem $(P)_\epsilon$ has a positive solution for $\epsilon>0$ small enough.   
\end{theorem}

\mbox{} \,\,\, In the proof of Theorem \ref{T3}, we follow the ideas found in \cite{OS}, because in that paper, problem $(P)_\epsilon$ has been considered with   $V$ verifying the same conditions explored in  \cite{Pino} and $f$ satisfying the above assumptions.

\end{document}